\documentclass[12pt]{article}
\usepackage{amsmath,amssymb,amsfonts,amsthm}
\usepackage[margin=1in]{geometry}
\usepackage{verbatim}

\usepackage{graphicx}
\usepackage{subcaption}

\theoremstyle{plain}

\title{Revisiting Mayer: Symmetric solutions for sporadic cases of the Map Color Theorem}
\author{Timothy Sun\\Columbia University}
\date{}

\begin{document}

\maketitle

\begin{abstract}
The original proof of the genus of the complete graphs $K_n$ depended on Mayer's \emph{ad hoc} solutions for $n = 18, 20, 23$. Recently, an improved solution for $K_{20}$ was found by the author. The purpose of this note is to use the theory of current graphs to interpret the aforementioned result and to provide new embeddings of $K_{18}$ and $K_{23}$. 
\end{abstract}

\section{Introduction}

At the heart of the Map Color Theorem of Ringel, Youngs, and others~\cite{Ringel-MapColor} is that the complete graph $K_n$ embeds in the orientable surface of genus
$$I(n) = \left\lceil\frac{(n-3)(n-4)}{12}\right\rceil.$$
The general approach is to restrict the search to symmetric embeddings derived from special kinds of directed graphs known as current graphs. The existence of appropriate current graphs seems more likely as the graphs get larger, but for smaller cases, we sometimes have to resort to improvised solutions---specialized current graphs (e.g. for $K_{30}$ in the original proof), or worse yet, embeddings that have no symmetries and must be written out in full. 

The latter situation occurred for $n = 18, 20, 23$, and the original proof depended on the \emph{ad hoc} embeddings of Mayer~\cite{Mayer-Orientables}. These remarkable embeddings were found by using many local modifications like flipping edges and adding handles, but they are, as Ringel~\cite{Ringel-MapColor} lamented, unstructured and asymmetric. Furthermore, Mayer presented these embeddings without their derivations---it is unlikely that these details can be rediscovered.

We give current graph constructions for each of these three cases. These solutions satisfy desirable properties: they are succinct, easily verifiable, and unified with the rest of the proof of Ringel~\emph{et al.} The principle for finding these embeddings comes from ``lifting'' a known current graph to an equivalent current graph of greater index. This suggests the possibility of other higher-index solutions for embeddings of the same derived graph. We then explore this larger space of current graphs to find one that enables a short modification of the derived embedding into a genus embedding of the target complete graph. 

\section{Current graphs}

We assume familiarity with current graphs, especially \S9 of Ringel~\cite{Ringel-MapColor}. For background in topological graph theory, see Gross and Tucker~\cite{GrossTucker}. The Heffter-Edmonds principle states that cellular, orientable embeddings of a graph are in one-to-one correspondence with \emph{rotation systems}, where each vertex is assigned a cyclic permutation, or a \emph{rotation}, of its incident edge ends. The rotation defines a cyclic ordering (say, counterclockwise) with respect to some fixed orientation of the surface, from which the faces of the embedding can be traced out. In the case of simple graphs, one can simply replace edge ends with corresponding neighbors. 

An \emph{index $k$ current graph} is a pair $(\phi, \alpha)$, where $\phi: D \to S$ is a cellular $k$-face embedding of a digraph $D$ on an orientable surface $S$ and $\alpha: E(D) \to \Gamma$ is a labeling of each arc of $D$ with a \emph{current}, an element of a cyclic group $\Gamma$. A vertex satisfies \emph{Kirchhoff's current law} (KCL) if its incoming currents sum to 0. Some vertices which do not satisfy KCL are called \emph{vortices}, and are labeled with lowercase letters. A \emph{circuit} is a boundary walk of a face of $\phi$, and the \emph{log} of the circuit records the currents encountered along the walk in the following manner: if we traverse arc $e$ along its orientation, we write down $\alpha(e)$; otherwise, we write down $-\alpha(e)$; if we encounter a vortex, we record its label.

We label each circuit $[0], [1], \dotsc, [k{-}1]$. To generate the embedding from the logs of these circuits, for each element $\gamma \in \Gamma$ in the group, take the log of circuit $[\gamma \bmod{k}]$ and add $\gamma$ to each of its non-letter elements to obtain the rotation at vertex $\gamma$. The rotation around each lettered vertex is ``manufactured'' so that the entire embedding is triangular and orientable. To facilitate this process, we make use of ``Rule R*'' (see, e.g., Ringel~\cite[\S 2.3]{Ringel-MapColor}), which states that for every edge $(i,k)$ in a triangulation, if the rotation at vertex $i$ is of the form $\dots \,\, j \,\, k \,\, l \,\, \dots$, then the rotation at $k$ is of the form $\dots \,\, l \,\, i \,\, j \,\, \dots$, for some other vertices $j$ and $l$. In the case of $K_{20}$, one lettered vertex splits in two because we obtain a permutation consisting of two cycles. These two vertices have disjoint neighborhoods.

The current graphs in this paper are of index 2 or 3. Furthermore, they satisfy the following additional ``construction principles'':
\begin{enumerate}
\item[(C1)] Each vertex is of degree 3 or 1, except the vortices, which have degree $k$.
\item[(C2)] KCL is satisfied at each non-vortex of degree 3.
\item[(C3)] The current incident with a degree 1 vertex is of order $2$ or $3$ in $\Gamma$. 
\item[(C4)] Every vortex is incident with each circuit.
\item[(C5)] Each nonzero element $\gamma \in \Gamma$ appears exactly once in the log of each circuit, except the order 2 element, which can be absent.
\item[(C6)] If circuit $[a]$ traverses arc $e$ along its orientation and circuit $[b]$ traverses $e$ in the opposite direction, then $\alpha(e) \equiv b-a \pmod{k}$. 
\end{enumerate}

By (C3), the log records the order $2$ element twice consecutively, producing several doubled edges in the resulting embedding. We follow the convention of condensing the two instances in the log together (this prevents the current graph from violating (C5)), drawing the corresponding edge without orientation, and omitting its degree-1 endpoint. 

These construction principles guarantee a triangular embedding of a graph that is nearly complete (see Ringel~\cite{Ringel-MapColor}). By the Euler polyhedral equation, the genus of any triangular embedding is $$\frac{|E|-3|V|+6}{6}.$$

Finding a current graph and deriving an embedding as described above is often referred to as the \emph{regular} part of the problem. The \emph{additional adjacency} step performs local operations to transform the embedding into one of the desired complete graph. Ideally, the latter should be short, preserving as much of the structure of the former as possible. 

\section{$n = 18$}

Let $S_k$ denote the orientable surface of genus $k$. Our starting point is the ``orientable cascade'' of Jungerman and Ringel~\cite{JungermanRingel-Octa} for embedding $K_{18}-9K_{2}$ in $S_{16}$. Their example can be interpreted as an index 2 current graph by considering its canonical double cover. We exhibit a different current graph that allows us to add the nine missing edges with two handles. 

\begin{figure}[ht]
\centering
\includegraphics[scale=0.8]{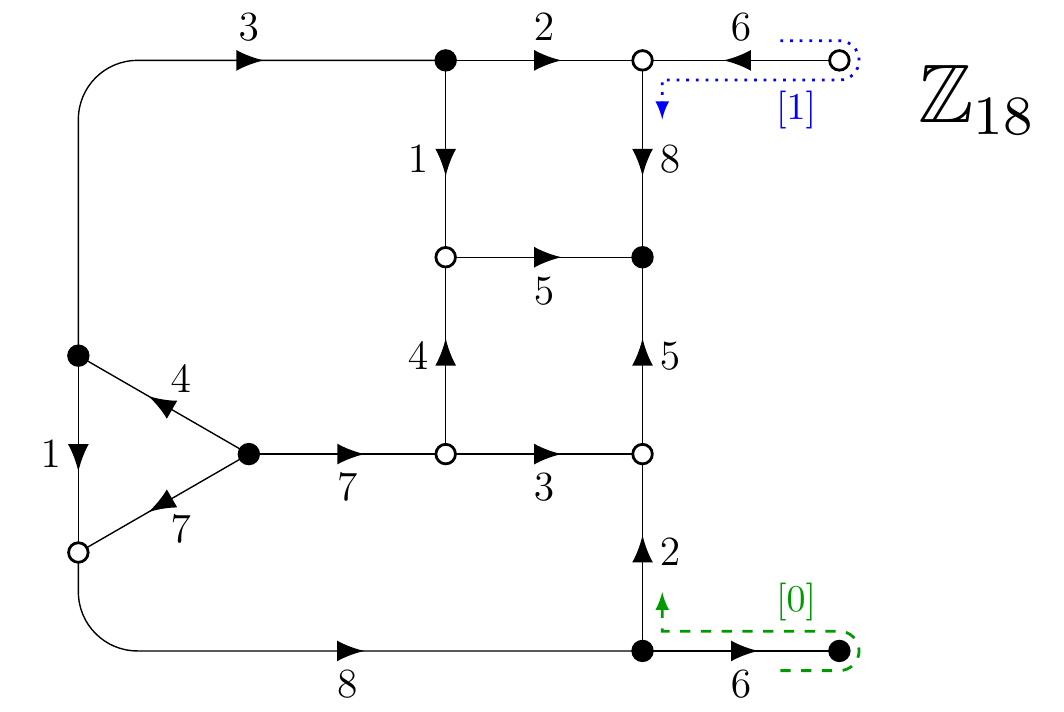}
\caption{An index 2 current graph for $K_{18}-9K_{2}$. The rotations at solid and hollow vertices are the edge ends arranged in clockwise or counterclockwise order, respectively.}
\label{fig-current-k18}
\end{figure}

\begin{figure}[ht]
\centering
\includegraphics[width=\textwidth]{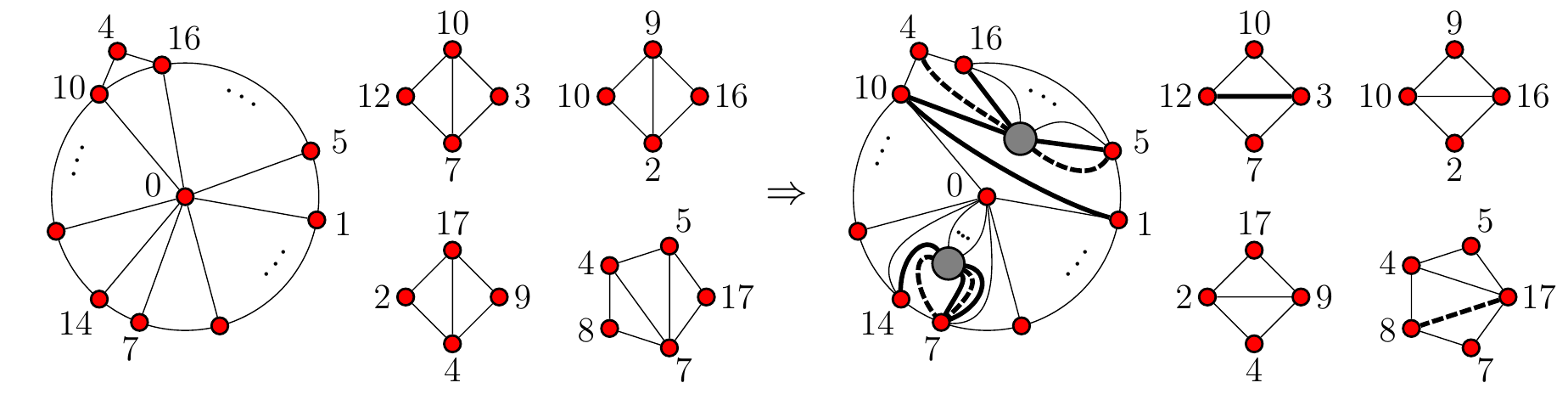}
\caption{Adding new edges by installing a handle, which is expressed as deleting the two gray disks and identifying their boundaries. The second handle adds only the solid bolded edges in this drawing.}
\label{fig-add-k18}
\end{figure}

The logs of the two faces in the current graph in Figure~\ref{fig-current-k18} are
$$\arraycolsep=3.6pt\begin{array}{rrrrrrrrrrrrrrrrr}
\lbrack0\rbrack. & 6 & 12 & 2 & 15 & 11 & 4 & 3 & 1 & 5 & 13 & 16 & 10 & 17 & 14 & 7 & 8  \\
\lbrack1\rbrack. & 12 & 6 & 8 & 13 & 14 & 3 & 5 & 10 & 16 & 15 & 1 & 11 & 7 & 4 & 17 & 2. \\
\end{array}$$
Edges of the form $(i, i+9)$ are missing from the graph. For all even $j$, the rotation system is of the form
\begin{equation*}
\begin{array}{rrrrrrrrrrrrrrrr}
j. & \dots & j{+}1 & j{+}5 & \dots & j{-}2 & j{-}8 & \dots & j{-}4 & j{+}7 & \dots \\
j{+}7. & \dots & j{+}3 & j{-}8 & j{-}6 & \dots. \\
\end{array}
\tag{\textasteriskcentered}
\label{eq-i2form}
\end{equation*}
We assert that this is enough to add the edges
$$(j{+}1, j{-}8), (j{+}3, j{-}6), (j{+}5, j{-}4), (j{+}7, j{-}2)$$
using an edge flip and one handle near vertex $j$. Figure~\ref{fig-add-k18} illustrates this operation for $j=0$, where some additional edge flips are included to add the missing edge $(8,17)$. After this augmentation, the substructure~(\textasteriskcentered) still holds for $j = 8$. Thus, we can add the remaining four edges 
$$(9, 0), (11, 2), (13, 4), (15, 6)$$
using the same construction, by incrementing all the numbers in Figure~\ref{fig-add-k18} by $8$. 

\section{$n = 20$}

\begin{figure}[ht]
\centering
\includegraphics[scale=.8]{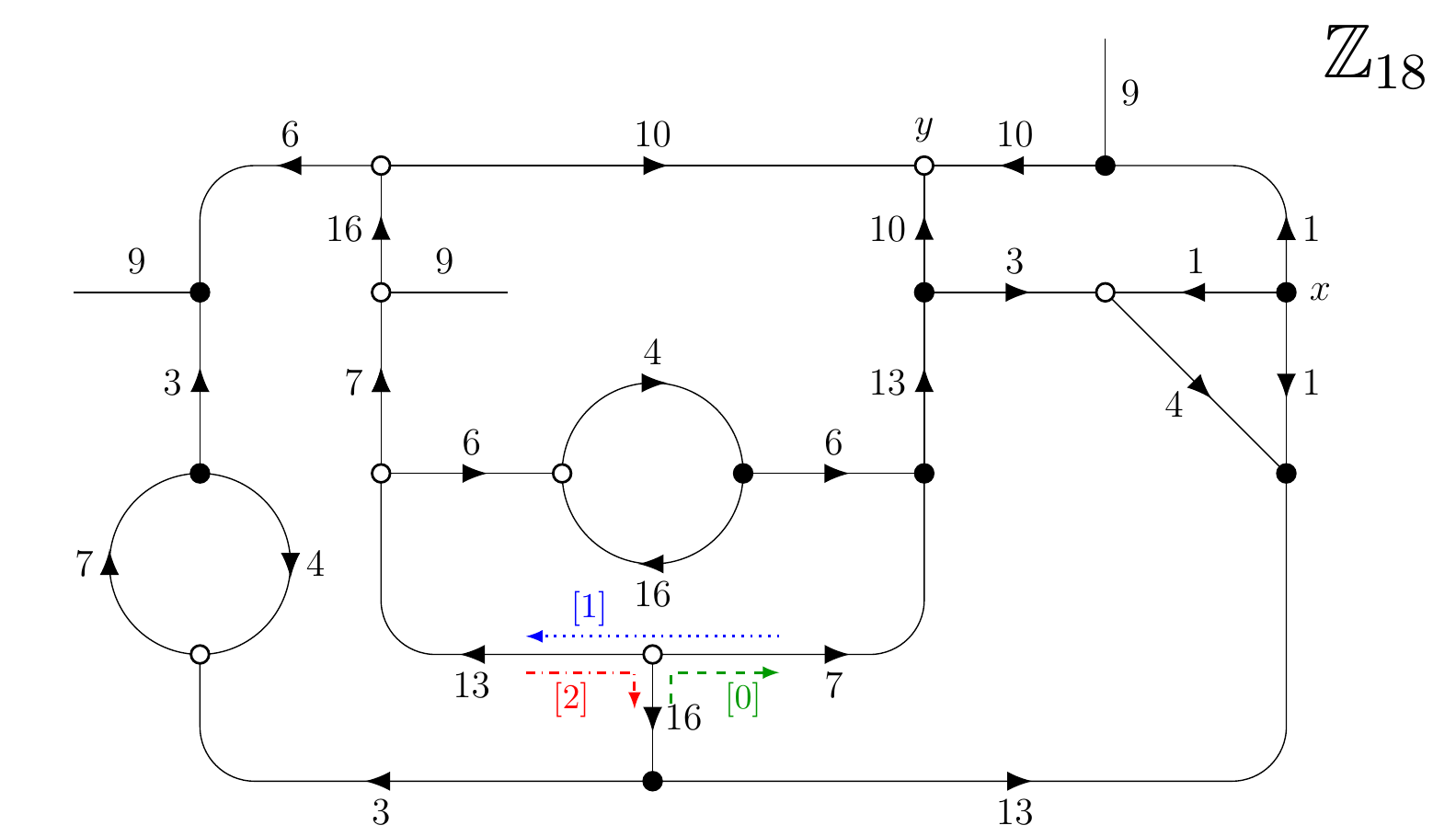}
\caption{An index 3 current graph for $G_{20}$.}
\label{fig-current-k20}
\end{figure}

\begin{figure}[ht]
\centering
\includegraphics[scale=.92]{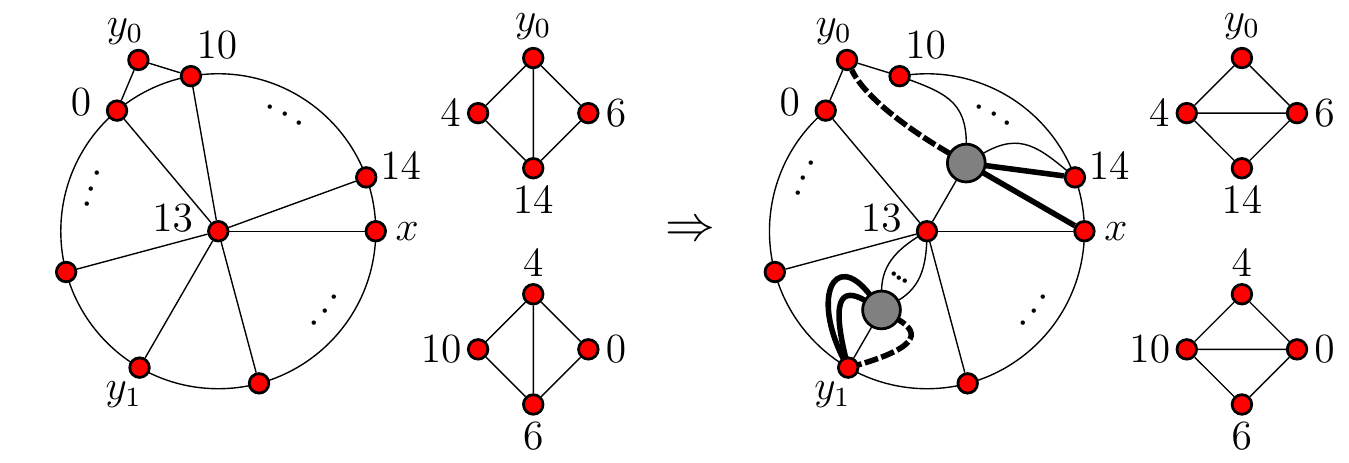}
\caption{The embedding is suitable for reusing Ringel and Youngs's construction.}
\label{fig-add-k20}
\end{figure}

The present author~\cite{Sun-FaceDist} gave a new solution for $K_{20}$ in $S_{23}$ that has the same properties as the index 1 current graphs of Ringel and Youngs (see \cite[\S7.5]{Ringel-MapColor}) for general $K_{12s+8}$. That solution can be derived from the current graph in Figure~\ref{fig-current-k20}, which yields
$$\arraycolsep=3.6pt\begin{array}{rrrrrrrrrrrrrrrrrrrrr}
\lbrack0\rbrack. & 2 & 7 & 12 & 16 & 4 & 6 & 13 & 10 & y & 8 & 9 & 17 & x & 1 & 5 & 3 & 14 & 11 & 15 \\
\lbrack1\rbrack. & 11 & 13 & 6 & 2 & 14 & 12 & 7 & 9 & 16 & 10 & y & 8 & 3 & 4 & 17 & x & 1 & 15 & 5 \\
\lbrack2\rbrack. & 5 & 16 & 13 & 14 & 17 & x & 1 & 10 & y & 8 & 6 & 15 & 4 & 7 & 3 & 9 & 12 & 2 & 11
\end{array}$$
as the logs of the three circuits. The embedded graph $G_{20}$ consists of vertices $$0, 1, \dotsc, 17, x, y_0, y_1,$$ where the numbered vertices form a clique, $x$ is adjacent to all the numbered vertices, and $y_0$ and $y_1$ are adjacent to the even and odd numbered vertices, respectively. One way of seeing why vortex $y$ splits into two vertices is to note that the currents entering $y$ are all even. The graph $G_{20}$ has been embedded in $S_{22}$, so the aim is to use one handle to connect the letters $x$, $y_0$, and $y_1$. We simply apply the additional adjacency step found in Ringel~\cite[\S7.6]{Ringel-MapColor}, which is shown in primal form in Figure~\ref{fig-add-k20}. After contracting the edge $(y_0,y_1)$, the embedding is now of $K_{20}$. 

\section{$n = 23$}

Ringel and Youngs (see \cite[\S7.2]{Ringel-MapColor}) found index 1 current graphs that generate triangular embeddings of $K_{12s+11}-K_5$ for $s \geq 1$. A special vortex forms a pattern in the rotation system that repeats every three vertices---this can be interpreted as an ordinary index 3 current graph. Unlike the previous case, our additional adjacency step differs from that of Ringel and Youngs for $s \geq 2$. It is a simpler one that resembles one of Jungerman and Ringel's~\cite{JungermanRingel-Minimal} steps for finding minimal triangulations on $12s{+}11$ vertices. 

\begin{figure}[ht]
\centering
\includegraphics[scale=0.8]{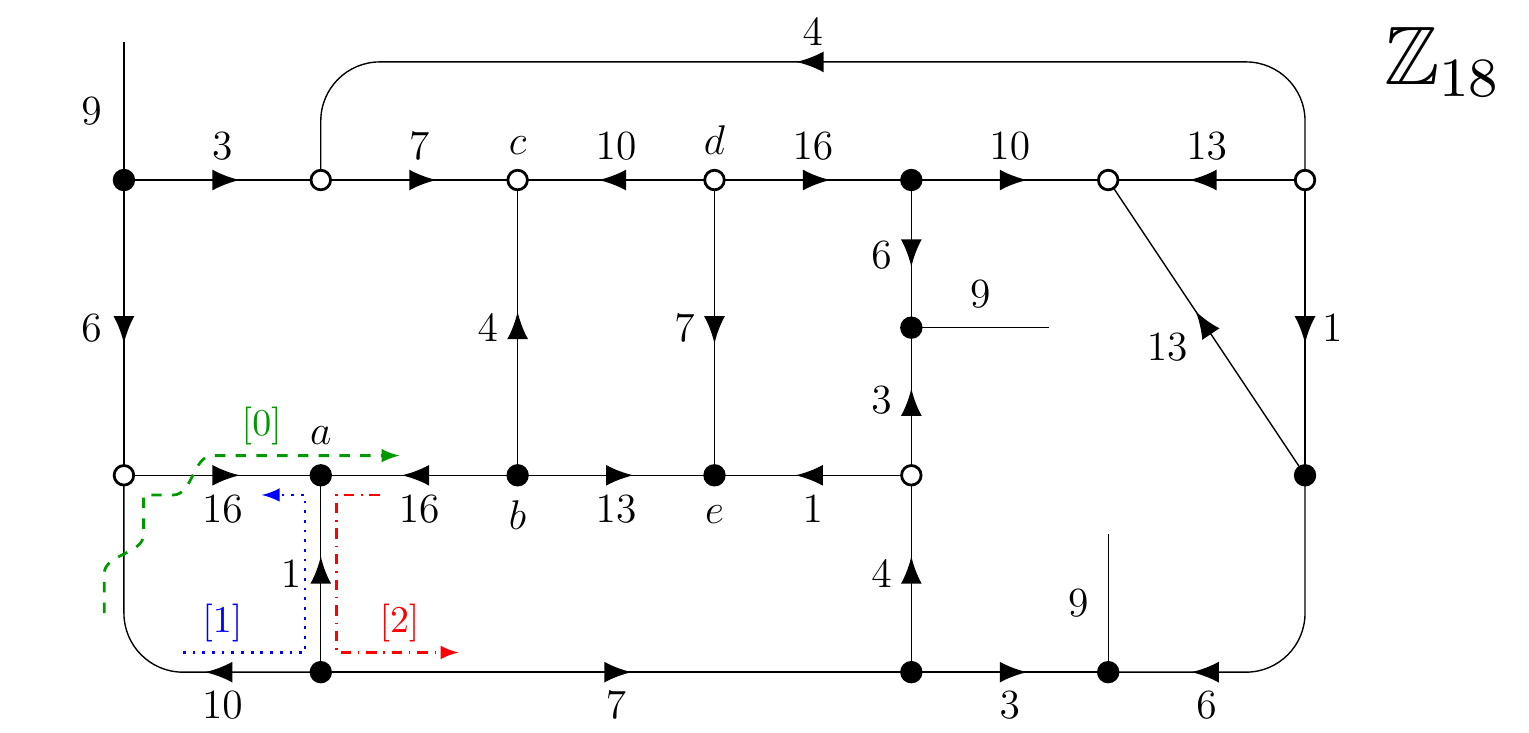}
\caption{An index 3 current graph for $K_{23}-K_5$.}
\label{fig-current-k23}
\end{figure}

\begin{figure}[ht]
    \begin{subfigure}[b]{0.49\textwidth}
        \centering
        \includegraphics[scale=0.92]{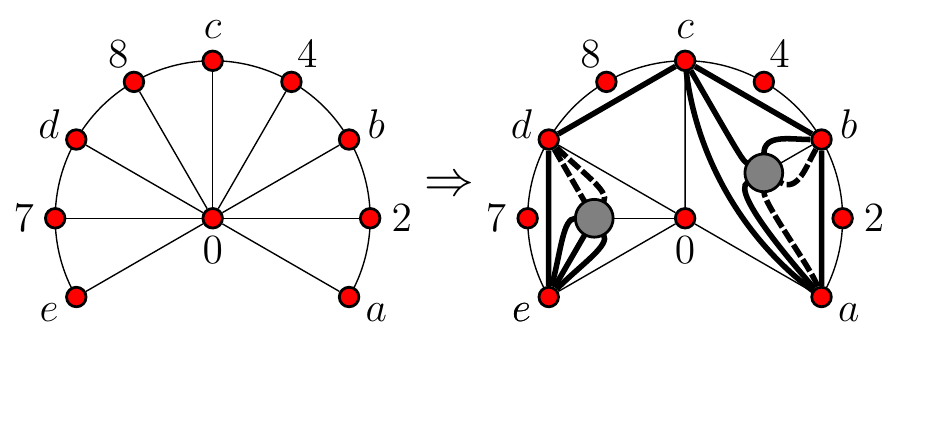}
        \caption{}
        \label{subfig-c1-a}
    \end{subfigure}
    \begin{subfigure}[b]{0.49\textwidth}
        \centering
        \includegraphics[scale=0.92]{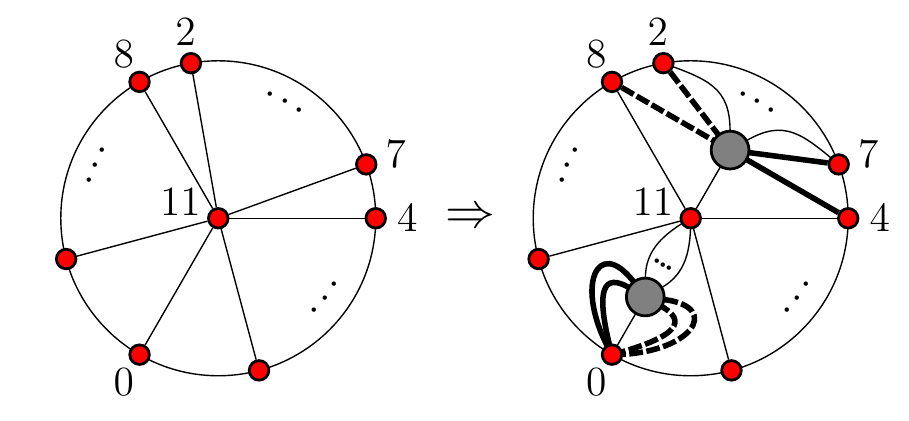}
        \caption{}
        \label{subfig-c1-b}
    \end{subfigure}
\caption{Connecting the letters (a) and then restoring the deleted edges (b).}
\label{fig-add-k23}
\end{figure}

The logs of the three circuits in Figure~\ref{fig-current-k23} are
$$\arraycolsep=3.6pt\begin{array}{rrrrrrrrrrrrrrrrrrrrrrr}
\lbrack0\rbrack. & 10 & 16 & a & 2 & b & 4 & c & 8 & d & 7 & e & 17 & 14 & 3 & 9 & 12 & 13 & 5 & 1 & 6 & 15 & 11 \\
\lbrack1\rbrack. & 8 & 1 & a & 2 & 12 & 9 & 3 & 7 & c & 14 & b & 13 & e & 11 & d & 16 & 10 & 5 & 17 & 4 & 15 & 6 \\
\lbrack2\rbrack. & 16 & a & 17 & 7 & 4 & 3 & 12 & 2 & d & 10 & c & 11 & 14 & 13 & 8 & 6 & 9 & 15 & 1 & e & 5 & b
\end{array}$$
which generate an embedding in $S_{30}$, where the missing edges are between lettered vertices. With one handle near vertex 0, we can add all ten of these edges at the cost of four edges incident with 0, as in Figure~\ref{fig-add-k23}(a). Figure~\ref{fig-add-k23}(b) shows how these four edges can be recovered using another handle near vertex 11.

\bibliographystyle{alpha}
\bibliography{../biblio}

\end{document}